\documentclass[11pt,a4paper]{article} 

\usepackage{fullpage}
\usepackage{amssymb,amsmath,stmaryrd,amsthm,nicefrac}
\usepackage{accents}

\usepackage{enumerate}
\usepackage[mathscr]{eucal}
\usepackage{epsfig}
\usepackage{pst-node,pst-text,pst-3d}
\usepackage{pstricks} 
\usepackage{pifont}

\usepackage{textcomp}   
\usepackage{multicol}

\theoremstyle{definition}
\newtheorem{theorem}{Theorem}
\newtheorem{lemma}{Lemma}
\newtheorem{remark}{Remark}

\newtheorem{example}{Example}

\def\1{\mathbf{1}}

\def\cB {\mathscr{B}}

\def\cG{\mathscr{G}}

\def\QQ{{\mathbb Q}}
\def\RR{{\mathbb R}}
\def\ZZ{{\mathbb Z}}
\def\NN{{\mathbb N}}
\def\CC{{\mathbb C}}

\def\ga{\mathfrak{a}}

\def \0{{\mathbf{0}}}

\def \im {\mathrm{i}}
\def \gcd{\mathrm{gcd}}

\begin{document}

\title{Polynomial representation of TU-games}
\author{Ulrich Faigle${}^{(a)}$ and Michel Grabisch${}^{(b)}$\\
$(a)$ Mathematisches Institut, Universit\"at zu K\"oln\\
  Weyertal 80, 50931 K\"oln, Germany\\
  {\tt faigle@zpr.uni-koeln.de}\\
$(b)$ Universit\'e Paris I Panth\'eon Sorbonne, Centre
  d'Economie de la Sorbonne\\ 106-112 Bd de l'H\^opital, 75013 Paris, France \\
    {\tt michel.grabisch@univ-paris1.fr}}

\maketitle

\begin{abstract}We propose in this paper a polynomial representation of TU-games,
  fuzzy measures, capacities, and more generally set functions. Our
  representation needs a countably infinite set of players and the natural
  ordering of finite sets of $\NN$, defined recursively. For a given basis of
  the vector space of games, we associate to each game $v$ a formal polynomial of
  degree at most $2^n-1$ whose coefficients are the coordinates of $v$ in the
  given basis. By the fundamental theorem of algebra, $v$ can be represented by
  the roots of the polynomial. We present some new families of games stemming
  from  this polynomial context, like the irreducible games, the multiplicative games and the cyclotomic
  games.
\end{abstract}

\section{Introduction}
The seminal notion of {\it fuzzy measure} introduced by Michio Sugeno
\cite{sug74} has led to many developments and applications, and has deep roots
in decision theory, game theory, as well as in discrete mathematics, operations
research and combinatorics in particular \cite{gra16}. Based on fuzzy measures,
Michio Sugeno has proposed a new type of integral, now called {\it Sugeno
  integral}. It is well known that fuzzy measures correspond to Choquet
capacities \cite{cho53}, and to the characteristic functions of monotone games
with transferable utilities (see, e.g., \cite{pesu03}). Particular fuzzy
measures can be found also in combinatorial optimization, as the rank function
of a polymatroid \cite{fuj05b}.

In this paper, we focus on the most general version of fuzzy measures, which
correspond to (the characteristic function of) TU-games. The traditional view of
TU-games (or games for short) is based on the fact that the set of games, for a
fixed set of players, say $N=\{1,\ldots,n\}$, is a vector space of dimension
$2^n-1$. This means that games can be added and multiplied by a scalar, and that
there exists a basis for this vector space. Moreover, game theorists have at
 their disposal the  whole arsenal of linear algebra to deal with games, defining linear
operators, computing their kernel, etc. The Shapley value \cite{sha53} is a
famous example of linear operator, mapping the set of games $\cG(N)$ to $\RR^N$,
but there are other important examples, as the  invertible linear operators from
$\cG(N)$ to $\cG(N)$, like the M\"obius transform, the interaction transform,
and so forth \cite{fagr14}. Each of these transform  implicitly defines a vector
basis of $\cG(N)$.  Another fundamental feature of games defined on a set $N$ of
players is the natural structure of the set of subsets of $N$ (called
coalitions) as  the Boolean lattice  $(2^N,\cup, \cap)$, where the supremum and the
infimum correspond to the union and intersection of sets. This naturally leads
to the notion of supermodular and submodular games, called in game theory convex
and concave games, respectively. A notable particular case is the case of {\it
  additive} games, i.e., those which are both supermodular and submodular.

The present paper proposes another view of games, which, up to our knowledge, 
has not been taken so far:  we represent a game as a formal polynomial,
the coefficients of which are precisely the quantities specifying a given game relative to 
a given basis. By the fundamental theorem of algebra, every polynomial of degree
$n$ has $n$ not necessarily pairwise different complex roots, and two different
polynomials have different roots, unless they differ only by a multiplicative
coefficient. It follows that a game of $n$ players can be represented by a
polynomial of degree at most $2^n-1$, and therefore can be represented by its
(at most) $2^n-1$ roots. This alternative to the vector basis representation
opens another world. As will be explained below, the polynomial view suggests 
to consider an infinite number of players, and addition of games is no more the
main operation, but multiplication of games, and as a consequence factorization
of games. These concepts seem to be entirely new in game theory. Another
consequence, coming from the abandon of a fixed set of players, is that the
Boolean lattice structure of the coalitions has to be replaced by something else,
which is a total order on the set of all finite coalitions. This makes the
notion of an additive game or  a supermodular game   inoperative.  Instead, we propose
the notion of  a multiplicative game. Therefore, the classical family of games  considered so far 
(convex games, simple games, monotone games, etc.) should be extended with
families stemming from the  polynomial context, for example the  family of 
cyclotomic polynomials, whose roots are the primitive roots of unity. Much
remains to be done in this unexplored region.

\section{The framework}
\subsection{Games and their bases}
We consider a countably infinite set of players $U$, and a bijection $\varphi$
from $U$ to $\NN$ (we assume throughout the paper that $0\not\in \NN$). We
denote by $[n]:=\{1,\ldots,n\}$ the set of $n$ first players, according to the
ordering $\varphi$.

For a given $n\in\NN$, classically a {\it TU-game} (or {\it game} for short) on
$[n]$ is a set function $v:2^{[n]}\rightarrow \RR$ such that $v(\varnothing)=0$.
In  our present model, we depart a little from this traditional view in two
respects:
\begin{enumerate}
\item For the sake of generality, we allow $v(\varnothing)\neq 0$  and thus 
  include general set functions,  which is also motivated by the Dempster-Shafer framework  \cite{sha76,gra16} in decision making. There, belief and
  plausibility functions, which are particular games, take nonzero values on the empty set.
\item We  typically take the range of $v$ to be $\QQ$ instead of $\RR$, which fully
  makes sense in most applications. 
\end{enumerate}
To summarize, a game $v$ on $[n]$ is a set function $2^{[n]}\rightarrow\QQ$.
The set of games on
$[n]$ is denoted by $\cG(n)$. 
It is a $2^n$-dimensional vector space, for which several canonical bases can be
defined. In particular:
\begin{enumerate}
\item The basis of {\it Dirac games} $(\delta_S)_{S\in
  2^{[n]}}$, with $\delta_S(T)=1$ iff $T=S$, and 0
  otherwise. The representation of a game $v$ in this basis is:
  \begin{equation}
v = \sum_{S\subseteq [n]}v(S)\delta_S.
  \end{equation}
  Hence, the coordinates of $v$ in this basis are simply $v(S)$, $S\subseteq
  [n]$.
\item The basis of {\it unanimity games} $(u_S)_{S\in
  2^{[n]}}$, with
  \begin{equation}
    u_S(T) = \begin{cases}
      1, & \text{ if } T\supseteq S\\
      0, & \text{ otherwise}.
    \end{cases}
  \end{equation}
  The representation of $v$ in this basis is
  \begin{equation}
v = \sum_{S\subseteq [n]}m^v(S)u_S
  \end{equation}
  The coefficients $m^v(S)$, $S\subseteq [n]$, correspond to
  the {\it M\"obius transform} \cite{rot64} of $v$:
  \[
m:\cG(n)\rightarrow\cG(n), \quad v\mapsto m^v
\]
with
\[
m^v(S) = \sum_{T\subseteq S}(-1)^{|S\setminus T|}v(T), \quad S\in2^{[n]}.
\]
\end{enumerate}
More generally, to each basis corresponds a transform, e.g., the interaction
transform, the Fourier transform, see \cite{fagr14,gra16} for details. In
general, we define transforms as linear invertible (i.e., bijective)
mappings $\Psi:\cG(n)\rightarrow \cG(n)$.  
We take
the convention to name the basis by the name of the corresponding transform: the
M\"obius basis, the interaction basis, etc. Note that under this convention the
basis of Dirac games is the identity basis.

We use the generic notation of $\cB^\Psi(n)=(b^\Psi_S)_{S\in
  2^{[n]}}$ for a basis corresponding to some transform
$\Psi$, hence with coefficients $\Psi^v(S)$:
\[
v = \sum_{S\subseteq [n]}\Psi^v(S)b^\Psi_S.
\]

Define $\cG:=\bigcup_{n\in \NN}\cG(n)$ (note that any game in $\cG$ has a finite
number of players), and $\cB^\Psi:=\bigcup_{n\in
  \NN}\cB^\Psi(n)$.
  
\subsection{Cardinality polynomials}\label{sec:cardinality-polynomials}
For the cooperative game $v\in \cG(n)$, define the $(n+1)$-dimensional parameter vector $\overline{v}$ of averages 
$$
   \overline{v}_k= \frac{1}{\binom{n}{k}}\sum_{S\subseteq [n], |S|=k} v(S)  \quad(k=0,1,\ldots,n) 
$$
and associate with it the formal polynomial 
$$
    \chi^v(x) =\sum_{k=0}^n \overline{v}_kx^k
$$
of degree $\leq n$. $v$ is reconstructible from $\chi^v(x)$ if and only if $v$ is a \emph{cardinality game}, i.e., 
$$
      v(S) = v(T) \quad\mbox{whenever $|S|=|T|$,} 
$$
and hence $v(S) = \overline{v}_{|S|}$ holds 
for all $S\subseteq [n]$. On the other hand, every polynomial can be interpreted as the cardinality polynomial of some cooperative game. 

We will now establish a polynomial representation of an arbitrary game $v$ from which $v$ can be reconstructed.

\subsection{The natural ordering of finite subsets of $\NN$}\label{sec:order}
We define a total order $\leqslant$ on $\bigcup_{n\in \NN}2^{[n]}$, the (countable) set of
all finite subsets of $\NN$: identify any finite subset $S\subseteq \NN$ with its
characteristic vector $\1^S\in\{0,1\}^\NN$, and define the bijection
\[
\eta:\bigcup_{n\in \NN}2^{[n]}\rightarrow \NN\cup\{0\}, \quad S\mapsto \sum_{k\in\NN}\1^S_k2^{k-1}
\]
Then $S\leqslant T$ iff $\eta(S)\leqslant \eta(T)$. This yields the order
\begin{multline}\label{eq:order}
  \varnothing,\{1\},\{2\},\{1,2\},\{3\},\{1,3\},\{2,3\},\{1,2,3\},\\ 
  \{4\},\{1,4\},\{2,4\},\{1,2,4\},\{3,4\},\{1,3,4\},\{2,3,4\},\{1,2,3,4\},\ldots
\end{multline}
Observe that this order is defined recursively: the sequence of subsets of
$[n+1]$ is the sequence of subsets of $[n]$ concatenated with the same sequence
where to each subset the element $n+1$ has been added. As it is based on the
order of natural numbers, we call it the {\it natural ordering}.

While the principle of this order is simple, it may be less obvious to find the
successor $S^+$ of a set $S$ in this order, when the rank of $S$ in the order is
getting large. Below is an algorithm to do it simply.
\begin{quote}
  \hrule
  
  {\sc Finding the successor $S^+$ of $S$}\\
  Define $n$ as the highest element in $S$\\
  {\bf if} $1\not\in S$ {\bf then} $S^+=S\cup\{1\}$\\
  \mbox{}$\quad$ {\bf else if} $2\not\in S$ {\bf then}
  $S^+=(S\setminus\{1\})\cup\{2\}$\\
  \mbox{}$\qquad$  {\bf else if} $3\not\in S$ {\bf then}
  $S^+=(S\setminus\{1,2\})\cup\{3\}$\\
  \mbox{}$\qquad \quad \cdots$\\
  \mbox{}$\qquad \qquad$  {\bf else if} $n-1\not\in S$ {\bf then}
  $S^+=(S\setminus\{1,\ldots,n-2\})\cup\{n-1\}$\\
  \mbox{}$\qquad\qquad\quad$ {\bf otherwise} $S^+=\{n+1\}$.

  \hrule
\end{quote}
This algorithm can be generalized for a shift of $2^k$ places in the ordered
sequence of sets ($k=0,1,\ldots$).
\begin{quote}
  \hrule
  
  {\sc Finding the $2^k$th successor $S^{2^k+}$ of $S$}\\
  Define $n$ as the highest element in $S$\\
  {\bf if} $k+1\not\in S$ {\bf then} $S^{2^k+}=S\cup\{k+1\}$\\
  \mbox{}$\quad$ {\bf else if} $k+2\not\in S$ {\bf then}
  $S^{2^k+}=(S\setminus\{k+1\})\cup\{k+2\}$\\
  \mbox{}$\qquad$  {\bf else if} $k+3\not\in S$ {\bf then}
  $S^{2^k+}=(S\setminus\{k+1,k+2\})\cup\{k+3\}$\\
  \mbox{}$\qquad \quad \cdots$\\
  \mbox{}$\qquad \qquad$  {\bf else if} $n-1\not\in S$ {\bf then}
  $S^{2^k+}=(S\setminus\{k+1,\ldots,n-2\})\cup\{n-1\}$\\
  \mbox{}$\qquad\qquad\quad$ {\bf otherwise} $S^{2^k+}=\{n+1\}$.

  \hrule
\end{quote}
Observe that by adequate combination of these two algorithms, one can determine
the $k$th successor of $S$, for any integer $k$.

Another interesting property is the following: take $S$, with integer
representation $\eta(S)$. Then $2\eta(S)$ corresponds to the set $S'$ where each
player $i\in S$ is replaced by $i+1$. For example $S=\{3,4\}$ has representation
$\eta(S)=12$, and 24 corresponds to the set $\{4,5\}$. More generally,
$2^k\eta(S)$ corresponds to a shift of $k$ positions for each player in $S$.

\medskip

Thanks to this ordering, a game in $\cG(n)$ can be represented by a vector in
$\QQ^{2^n}$, where coordinates are the coordinates of this game in some
basis. Specifically, with basis $\cB^\Psi$, the $k$th coordinate of $v$ is
$\Psi^v(\eta^{-1}(k))$, for any $k\in\NN\cup\{0\}$.
Conversely any vector in $\QQ^{2^n}$ corresponds to a unique game in
$\cG(n)$, provided the basis is chosen. Therefore, in the remainder of the
paper, when considering a game $v\in\cG(n)$ represented by $\Psi^v$, we may use either
the transform notation $\Psi^v(S)$, $S\in 2^{[n]}$, or the vector notation
$\Psi^v_k$, with $k=0,1,\ldots,2^n-1$, the correspondance being given by $k=\eta(S)$.

\subsection{Equivalence classes of games}
Consider a game $v\in\cG(n)$, with representation $\Psi^v\in\QQ^{2^n}$. We may
consider $\Psi^v$ as a vector in $\QQ^{\bigcup_{n\in\NN}2^{[n]}}$, by padding
  $\Psi^v$ with zeroes. Let us denote the vector with padding zeroes by
  $\overline{\Psi^v}$. 

Consider a given transform $\Psi$ and $\cB^\Psi$. We define the equivalence
relation $\equiv^\Psi$ on $\cG$ as follows:
\[
v\equiv^\Psi v' \Leftrightarrow \overline{\Psi^v}=\alpha\overline{\Psi^{v'}},
\]
for some $\alpha\in\RR\setminus\{0\}$. 
Denote by $[v]^\Psi$ the equivalence class of $v$. Considering an equivalence
class $[v]^\Psi$ with $v\in\cG(n)$, its {\it canonical representative} is the
game $\underline{v}\in[v]^\Psi$ with the least number of zeroes and with
$\Psi^{\underline{v}}(\underline{S})=1$, where $\underline{S}$ is the last set $S$
in the order (\ref{eq:order}) such that $\Psi^{\underline{v}}(S)\neq 0$ (games
satisfying the former property are called {\it minimal} and those satisfying the
latter are called {\it normalized}. Observe, however, that these notions depend
on the chosen basis). Observe that
$\underline{v}\in\cG(n')$, with $n'\leqslant n$.  We denote by
$\cG/_{\equiv^\Psi}$ the corresponding quotient set of games.

\begin{example}\label{ex:1}
Consider the M\"obius basis, $n=3$ and the game $v=u_{\{1,2\}}+2u_{\{3\}}$ in
$\cG(3)$, which has representation $m^v=(0 \ 0 \ 0 \ 1 \ 2 \ 0 \ 0 \ 0)$ in the
M\"obius basis. Then $v\equiv^m v'$ with $v'\in\cG(4)$ having representation
\[
m^{v'}=(0 \ 0 \ 0 \ \alpha \ 2\alpha \ 0 \ 0 \ 0\ 0 \ 0 \ 0 \ 0 \ 0 \ 0 \ 0 \ 0)
\]
with $\alpha\neq 0$, 
and similarly with $v'\in\cG(n')$, for any $n'>3$, padding $m^v$ with $2^{n'}-2^3$
zeroes.

\hfill$\Diamond$

\end{example}

\medskip

\subsection{Polynomial and algebraic representations of games}
Consider a game $v\in\cG(n)$, with representation $\Psi^v$ in basis
$\cB^\Psi$. We consider the formal polynomial $P^{v,\Psi}$ in the indeterminate $x$ with
rational coefficients:
\[
P^{v,\Psi}(x) = \sum_{i=0}^{2^n-1}\Psi^v_ix^{i}
\]
and call it the {\it polynomial representation} of $v$.
Let us denote by $d$ the degree of this polynomial, with $d\leqslant 2^n-1$.
We know by the fundamental theorem of algebra that any such polynomial has $d$
roots $r_1,\ldots,r_d$ in $\CC$,  possibly 
non distinct, which are algebraic numbers (see, e.g., \cite{podi98}). We call
the multiset $\{r_1,\ldots,r_d\}$ the {\it algebraic representation} of $v$,
denoted by $\ga(v,\Psi)$. Observe
that it depends on the chosen basis.
Conversely, any family of $d$
algebraic numbers $r_1,\ldots,r_d$ defines a unique polynomial of degree
$d$ with rational coefficients $a_0,a_1,\ldots,a_{d-1},a_d =1$
\[
P(x)=(x-r_1)\cdots(x-r_d)= a_0+a_1x+\cdots + a_{d-1}x^{d-1}+x^{d}
\]
This polynomial is {\it monic}, i.e., the coefficient of its highest degree
monomial is 1. Summarizing, we have established the following fact.
\begin{lemma}
Given a transform $\Psi$, each game $v$ in $\cG$ has a (unique) algebraic
representation $\ga(v,\Psi)$. The algebraic representation
corresponds to a unique element of $\cG/_{\equiv^\Psi}$, i.e., all games
$v'\in[v]^\Psi$ have the same algebraic representation as $v$, and only those games. 
\end{lemma}
We denote by $\QQ[x]$ the set of all polynomials with rational coefficients, and
by $\QQ^1[x]$ the set of monic polynomials. When $P(x)\in\QQ[x]$, we also say
that $P(x)$ is a polynomial over $\QQ$. 

We have established a bijection
between $\QQ^1[x]$ and $\cG/_{\equiv^\Psi}$.

\begin{remark}
It is important to note that,  in contrast to most classical approaches, the
notion of permutation of players does not make sense here. In the classical
view, for a given game $v$, permuting the players by a permutation $\sigma$ is immaterial as it
suffices to consider the game $v^\sigma$ defined by
$v^\sigma(\sigma(S))=v(S)$. Then most of the classical notions like the Shapley
value are invariant to permutations. In our framework, recall that an order on
the players has been fixed by the bijection $\varphi$ from $U$ to
$\NN$. Switching players 1 and 2 for example will change the polynomial and the
algebraic representations of the two games, and they will {\it not} be in general
identical up to a permutation. Take for example the game with two players $v(1)=1$,
$v(2)=-1$, $v(12)=1$. Its polynomial and algebraic representations are:
\[
P^{v,Id}(x)=x -x^2+x^3, \quad \ga(v,Id) = \left\{0,\frac{1+\im\sqrt{3}}{2},\frac{1-\im\sqrt{3}}{2}\right\}
\]
while for the permuted game $v^\sigma$ we have:
\[
P^{v^\sigma,Id}(x) = -x+x^2+x^3, \quad \ga(v^\sigma,Id) = \left\{0,\frac{-1+\sqrt{5}}{2},\frac{-1-\sqrt{5}}{2}\right\},
\]
which are completely different.
\end{remark}

We recall the Vieta formula relating coefficients of a polynomial of degree $n$
to its roots:
\[
\sum_{1\leqslant i_1<i_2<\cdots<i_k\leqslant n}\Big(\prod_{j=1}^kr_{i_j}\Big)=
(-1)^k\frac{a_{n-k}}{a_n} \quad k=1,\ldots,n
\]
In particular, $\Pi_{i=1}^nr_i=(-1)^n\frac{a_0}{a_n}$, and $\sum_{i=1}^nr_i =
-\frac{a_{n-1}}{a_n}$. 
\begin{example}[Example~\ref{ex:1} ct'd]\label{ex:2}
  Consider $v\in\cG(3)$ of Ex.~\ref{ex:1}, again in the M\"obius basis. The associated polynomial is
  \[
P^{v,m}(x) = x^3 + 2x^4 = x^3(1+ 2x)
\]
whose roots are $-\nicefrac{1}{2}$ and 0 with multiplicity 3. Hence
\[
\ga(v,m) = \{-\nicefrac{1}{2}, 0, 0, 0\}.
\]
The same game in the identity basis has coordinates
$(0,\ 0,\ 0,\ 1,\ 2,\ 2,\ 2,\ 3)$,  which yields the polynomial
\[
P^{v,v}(x) = x^3 + 2x^4 + 2 x^5 + 2 x^6 + 3x^7 = x^3(1+ 2x+2x^2+2x^3+3x^4),
\]
whose roots are 0 with multiplicity 3 and
\begin{align*}
r_4=-0.568455251719284-0.334436436804741\im, \quad r_5=-0.568455251719284+0.334436436804741\im\\
r_6=0.23512191838595-0.843220258232994\im, \quad r_7=0.23512191838595+0.843220258232994\im.
\end{align*}
Exact expressions are available but very complex. The example shows the
importance of the basis for the simplicity of the representation.

\hfill$\Diamond$

\end{example}
Let us give other examples:
\begin{enumerate}
\item In the M\"obius basis, the unanimity game $u_S$ is represented by the
  monomial $x^{\eta(S)}$. In particular:
  \begin{equation}\label{eq:u}
P^{u_\varnothing}(x)=1, \ P^{u_{\{1\}}}(x)=x, \ P^{u_{\{2\}}}(x)=x^2,
\ldots, P^{u_{\{k\}}}(x)=x^{2^{k-1}}, \ldots
  \end{equation}
\item Consider the simple game $v$ whose minimal winning coalitions are $\{1\}$
  and $\{2\}$, i.e., $v(S)=1$ whenever $S\supseteq \{1\}$ or $S\supseteq \{2\}$, and $v(S)=0$
  otherwise. The representation of $v$ in the M\"obius basis is
  \[
v = u_{\{1\}} + u_{\{2\}} - u_{\{1,2\}},
\]
hence
\[
P^{v,m}(x) = x+x^2-x^3 = x(1+x-x^2).
\]
We get $\ga(v,m) = \{0,\frac{1+\sqrt{5}}{2},\frac{1-\sqrt{5}}{2}\}$. Thanks to
the property of the natural order, we can find immediately the polynomial and
algebraic representations of the game $v'$
\[
v' = u_{\{2\}}+u_{\{3\}} -u_{\{2,3\}}
\]
obtained from $v$ by shifting each player in the basis representation by one
position. We get
\[
P^{v',m}(x) = x^2+x^4-x^6=x^2(1+x^2-x^4),
\]
hence letting $y=x^2$, the roots of $P^{v',m}$ are the square roots of the roots
of $P^{v,m}$:
\[
\ga(v',m)=\left\{0,0,\sqrt{\frac{1+\sqrt{5}}{2}},-\sqrt{\frac{1+\sqrt{5}}{2}},\im\sqrt{\frac{-1+\sqrt{5}}{2}},-\im\sqrt{\frac{-1+\sqrt{5}}{2}}\right\}.
\]
Similarly, one can find the polynomial and algebraic representations of all
games obtained by a shift of $k$ positions in the basis representation of $v$.
\end{enumerate}

\section{Factorization and irreducible games}
  Traditional TU game theory emphasizes the role
of linear mappings, taking advantage of the structure of vector space of the set
of games of $\cG(n)$, for a fixed $n$: in addition to  the M\"obius transform cited above,
the Shapley value is another well-known example of a linear mapping, and all questions
around it, like finding the set of games with same Shapley value, reduce to
classical questions of linear algebra \cite{fagr14}.

The polynomial view of games we adopt in this paper orients us towards a
different world. First, the universe to be considered is no more $\cG(n)$ for a
given $n$, but the whole set $\cG$ of all games on finite sets. The second essential
difference is that polynomials form a commutative ring, hence the new operation
here is {\it multiplication}. As far as we know, multiplication of games has
never been considered so far. This is what the polynomial view offers. 
We start by recalling the necessary background on the factorization of
polynomials.

\subsection{Factorization of polynomials}
A polynomial $P(x)\in\QQ[x]$ is {\it irreducible} or {\it prime} if
it cannot be factored into a product of polynomials
$P_1(x),\ldots,P_k(x)\in\QQ[x]$, i.e., $P(x)=P_1(x)\cdots P_k(x)$. If
$P(x)=F(x)G(x)$ with $F(x),G(x)\in\QQ[x]$, then $F(x)$ {\it divides} $P(x)$. 

We give below classical results on irreducible polynomials (see, e.g.,
\cite{podi98}). A {\it unit} is a constant polynomial (i.e., of degree
0). Polynomials are {\it relatively prime} if they have only units as common
factors.
\begin{theorem}
Let $F(x),G(x)\in\QQ[x]$. If $P(x)\in\QQ[x]$ is irreducible and divides
$F(x)G(x)$, then $P(x)$ divides either $F(x)$ or $G(x)$. 
\end{theorem}
\begin{theorem}[uniqueness of factorization]\label{th:2}
Any polynomial $P(x)=a_0+a_1x+\cdots+a^nx^n$ in $\QQ[x]$ of degree at least one
can be factored into a product
\[
P(x) = a_nP_1(x)\cdots P_r(x)
\]
where the $P_i(x)$ are irreducible monic polynomials over  $\QQ$. Moreover, the
factorization is unique up to the order. 
\end{theorem}
With $p,q\in \ZZ$, $p$ divides $q$ is denoted by $p\mid q$. If $p$ does not divide
$q$, we write $p\nmid q$.  
\begin{theorem}[Eisenstein irreducibility criterion]\label{th:3}
Let $p$ be a prime number and $f(x)=a_0+a_1x+\cdots+a^nx^n$ a polynomial over
$\ZZ$ such that $p^2\nmid a_0$, $p\nmid a_n$ and $p\mid a_i$ for
$i=0,1,\ldots,n-1$. Then $f(x)$ is irreducible over $\QQ$. 
\end{theorem}

There is a fairly large amount of research about factorization of polynomials, going
back to Kronecker, and yielding practical algorithms able to factorize
polynomials of degree about 1000 (algorithms by Kronecker, Zassenhaus,
Lenstra-Lenstra-Lov\'asz, etc.). As shown by Gauss, factorization in $\QQ[x]$
amounts to factoring in $\ZZ[x]$. More precisely, let us consider a polynomial
$P(x)$ in $\QQ[x]$. It can be written in the following form:
\[
P(x) = \frac{Q(x)}{c}
\]
with $Q(x)\in\ZZ[x]$, $c\in\ZZ$, so that $\pm 1$ is the greatest common divisor
of all coefficients in $Q(x)$ (such a polynomial is called {\it
  primitive}). Then $Q(x)$ can be factored over $\QQ$ iff it can be factored
over $\ZZ$. It follows that the factorization problem over $\QQ$ amounts to the
factorization problem over $\ZZ$. 

\subsection{Products of games}
Let $v,w\in\cG$ be two games whose polynomial representation w.r.t. the
transform $\Psi$ are two monic polynomials $P^{v,\Psi}(x),P^{w,\Psi}(x)$. The {\it
  product} of $v$ and $w$ w.r.t. $\Psi$, denoted by $v\times^\Psi w$, is the
canonical representative in $\cG$ of the polynomial $P^{v,\Psi}(x)Q^{w,\Psi}(x)$.

By extension, consider two arbitrary games $v,w\in\cG$ whose polynomial
representations $P^{v,\Psi}(x),P^{w,\Psi}(x)$ are of degree $n,m$ respectively, with
coefficients $a_n,b_m$. Let $z$ be the canonical representative in $\cG$ of the
polynomial  $P^{v,\Psi}(x)Q^{w,\Psi}(x)$. Then the product of $v,w$ is defined as
follows:
\[
v\times^\Psi w = a_nb_mz.
\]
Note that this amounts to taking the game with representation
$P^{v,\Psi}(x)Q^{w,\Psi}(x)$ having the less padding zeroes (i.e., being minimal).

\begin{example}\label{ex:3}
Consider the two following games $v\in\cG(3)$, $w\in\cG(2)$ in the M\"obius
basis:
\begin{center}
  \begin{tabular}{|c|cccccccc|}\hline
    $S$ & $\varnothing$ & 1 & 2 & 12 & 3 & 13 & 23 & 123 \\ \hline
    $m^v(S)$ & 0 & 0 & 0 & $-1$ & 0 & 2 & $-1$ & 0 \\ \hline 
  \end{tabular}
  \hfill
  \begin{tabular}{|c|cccc|}\hline
    $S$ & $\varnothing$ & 1 & 2 & 12 \\ \hline
    $m^w(S)$ & 1 & 0 & -1 & 1\\ \hline
  \end{tabular}
\end{center}
We obtain:
\begin{align*}
  P^{v,m}(x) & = -x^3 + 2x^5 - x^6\\
  P^{w,m}(x) & = 1-x^2+x^3\\
  P^{v,m}(x)P^{w,m}(x) & = -x^3 + 3x^5 - 2x^6 - 2x^7 + 3x^8 - x^9.
\end{align*}
Therefore, the product $v\times^mw$ is a game in $\cG(4)$ with M\"obius
representation
\begin{center}
  \begin{tabular}{|c|cccccccccc|}\hline
    $S$ & $\varnothing$ & 1 & 2 & 12 & 3 & 13 & 23 & 123 & 4 & 14 \\ \hline
    $m^{v\times^mw}$ & 0 & 0 & 0 & $-1$ & 0 & 3 & $-2$ & $-2$ & 3 & $-1$\\ \hline
  \end{tabular}
\end{center}

\hfill$\Diamond$

\end{example}
The above example shows that taking the product of two games with $n$ and $m$
players may yield a new game with more players. Since the highest possible
degrees of the polynomials representing $v$ and $w$ are $2^n-1$ and $2^m-1$
respectively, the highest possible degree of the product is $2^n+2^m-2$, and
supposing $n\leqslant m$, this is smaller than $2^{m+1}$. Consequently, there is
at most one new player $m+1$, but there could be no new player at all (take,
e.g., $v$ represented by $x^2$ ($n=2$) and $w$ represented by $x^8$ ($m=4$)).

While this fact may appear as counterintuitive in the traditional framework, 
an interpretation is possible as follows. Let us take a game $v\in\cG(n)$ and a
basis $\cB^\Psi$ for its representation, and let us assume that $n$ minimizes
the number of padding zeroes for $v$, i.e., $v$ is minimal. We may say that all
players in $[n]$ are active, since each player is involved in at least one
coalition with a nonzero coordinate in the basis. Then, the other players
$n+1,n+2,\ldots$ may be considered as inactive, though present. Let us now take
two games $v$ and $w$ with $n$ and $m$ active players respectively, and suppose
$n\leqslant m$. This implies $[n]\subseteq [m]$, so that the set of active
players is $[m]$. Multiplying $v$ and $w$ has the effect of possibly triggering
a new player $m+1$ into becoming active.

Taking two games $v,w$ with polynomial representation $P^{v,\Psi}(x),
P^{w,\Psi}(x)$, what is the effect of multiplying the monomial $a_kx^k$ of $v$
with the monomial $b_lx^l$ of $w$? Indices $k,l$ refer to coalitions
$S=\eta^{-1}(k)$ and $T=\eta^{-1}(l)$, respectively. The product of the two
monomials yields $a_kb_lx^{k+l}$, which represents a shift of $l$ positions from
$S$ in the natural ordering, or equivalently, $k$ positions from $T$. The two
algorithms given in Section~\ref{sec:order} permit to find the $l$th successor of
$S$, or the $k$th successor of $T$. A simple application of these algorithms
yields to the following interesting fact:
  \[ \text{If }S\cap T=\varnothing,\text{ then the
    multiplication yields }S\cup T\]

\subsection{Factorization of games}
We say that a game $v\in\cG$ is {\it irreducible} or {\it prime} w.r.t. the
transform $\Psi$ if the corresponding polynomial $P^{v,\Psi}(x)$ is irreducible.
Theorem~\ref{th:2} is fundamental and permits to obtain a unique factorization
of games. We restate it in the framework of games.
\begin{theorem}\label{th:th4}
  Let $v\in\cG$ be a game and $\Psi$ be a transform. Then $v$ can be written as
  a product of irreducible games w.r.t. $\Psi$ which are minimal and normalized
  (hence canonical representative of some classes):
  \[
v = a_n\big(\underline{v_1}\times^\Psi
\underline{v_2}\times^\Psi\cdots\times^\Psi \underline{v_r}\big)
\]
Moreover, the factorization is unique up to the order.
\end{theorem}
It follows that the algebraic representation of $v$ is:
\[
\ga(v,\Psi) =
\ga(\underline{v_1},\Psi)\sqcup\ga(\underline{v_2},\Psi)\sqcup\cdots\sqcup\ga(\underline{v_r},\Psi)  
\]
where $\sqcup$ indicates the union (concatenation) of multisets.

The Eisenstein criterion (Theorem~\ref{th:3}) gives a sufficient condition for
irreducibility of games. Standard algorithms for factorization of polynomials
can be used to factorize games. 

\begin{example}[Ex.~\ref{ex:1} cont'd]
Let us consider the game $v$ in Example~\ref{ex:2}. When represented
in the M\"obius basis, its factorization was already obtained above: $P^{v,m}(x)
=  x^3(1+2x)$. In the identity basis, the polynomial is
\[
P^{v,v}(x) = x^3(1 + 2x + 2x^2 + 2x^3 + 3x^4).
\]
It can be checked that $1+2x+2x^2+2x^3+3x^4$ is irreducible. 

\hfill$\Diamond$

\end{example}
\begin{example}
  Let us consider $v\in\cG(5)$ defined by
  \[
v= u_{1}+u_{12} - 2u_{3} +2u_{13}- u_{24} + u_{35} -2u_{245} +u_{2345}
  \]
  Its polynomial representation in the M\"obius basis is
  \[
P^{v,m}(x) = x +x^3 -2x^4 + 2x^5 + x^{10} + x^{20} -x^{26} - x^{30}.
\]
Observe that the Eisenstein criterion does not permit to decide whether it is
irreducible or not. Its factorization in irreducible polynomials is
\begin{multline*}
P^{v,m}(x) =
-x(x^2-x+1)(x^{27}+x^{26}-x^{24}+x^{22}+x^{21}-x^{19}-x^{18}-x^{17}+x^{15}+x^{14}-x^{12}\\
-x^{11}+x^{9}+x^{8}-x^{7}-2x^6-x^5+x^4+2x^3-x^2-x-1).
\end{multline*}
Observe that $v$ is factored in 3 games
$\underline{v_1},\underline{v_2},\underline{v_3}$, where the first one is the
one-player game $u_{1}$, the second one is the 2-player game
$u_{\varnothing}-u_{1}+u_{2}$, and the third one is a game with 5 players. 

\hfill$\Diamond$

\end{example}

\section{Multiplicative games}
Traditional TU game theory emphasizes the role of {\it additive games}, i.e., games
satisfying the property
\[
v(S\cup T) = v(S) + v(T)
\]
whenever $S,T$ are disjoint subsets. They are characterized by the fact that
their M\"obius representation has nonzero coefficients only for singletons,
which leads to the obvious generalization of {\it $k$-additive games}, for which
the M\"obius representation has nonzero coefficients only for subsets of at most
$k$ elements \cite{gra96f}.

Additive and $k$-additive games do not seem to exhibit interesting properties in
the polynomial representation, apart that their maximum degree is known in the
M\"obius basis: $\binom{n}{1} + \binom{n}{2} + \cdots + \binom{n}{k}$ for a
$k$-additive game. In particular, there does not seem to be any way in general
to find the roots of an additive game.

\medskip

As we have seen, the polynomial representation is multiplicative in nature, not
additive. Therefore, we introduce the notion of a multiplicative game. 
Recall that any additive game $v\in\cG(n)$ (assuming $v(\varnothing)=0$) can be written as
\[
v = \alpha_1u_{\{1\}} + \alpha_2u_{\{2\}} +\cdots+\alpha_nu_{\{n\}}
\]
with $\alpha_1,\ldots,\alpha_n\in\QQ$.
By analogy, we say that $v\in\cG(n)$ is a {\it multiplicative game} if it has the form
\begin{equation}
v =\beta
(u_{\{1\}}-\alpha_1u_\varnothing)\times^m(u_{\{2\}}-\alpha_2u_\varnothing)\times^m\cdots\times^m(u_{\{n\}}-\alpha_nu_\varnothing), 
\end{equation}
with $\beta,\alpha_1,\ldots,\alpha_n\in\QQ$.
Recalling (\ref{eq:u}), its polynomial representation in the M\"obius basis is immediate:
\[
P^{v,m}(x) = \beta(x-\alpha_1)(x^2-\alpha_2)\cdots(x^{2^{n-1}}-\alpha_n)
\]
whose roots are the $2^{k-1}$th roots of unity multiplied by
$(\alpha_k)^{\frac{1}{2^{k-1}}}$, for $k=1,\ldots,n$. This can indeed represent
a game in $\cG(n)$, as the degree of this polynomial is $2^n-1$.

A simple example permits to see the interest of multiplicative games. Let us
take $n=3$ and consider an arbitrary multiplicative game with polynomial
representation
\[
P^{v,m}(x) = (x-\alpha_1)(x^2-\alpha_2)(x^4-\alpha_3).
\]
Multiplying out, one obtains the form
\[
P^{v,m}(x) = -\alpha_1\alpha_2\alpha_3 + \alpha_2\alpha_3x + \alpha_1\alpha_3x^2
-\alpha_3x^3 +\alpha_1\alpha_2x^4- \alpha_2x^5-\alpha_1x^6+x^7.
\]
So the coordinates of $v$ in the M\"obius basis are
\begin{center}
  \begin{tabular}{|c|cccccccc|}\hline
    $S$ & $\varnothing$ & 1 & 2 & 12 & 3 & 13 & 23 & 123\\ \hline
    $m^v(S)$ & $-\alpha_1\alpha_2\alpha_3$ & $\alpha_2\alpha_3$ &
    $\alpha_1\alpha_3$ & $-\alpha_3$ & $\alpha_1\alpha_2$ & $-\alpha_2$ &
    $-\alpha_1$ & 1\\ \hline
  \end{tabular}
\end{center}
In summary, we have obtained
\[
m^v(S) = (-1)^{|S|+1}\prod_{i\not\in S}\alpha_i,
\]
with the convention $\prod_\varnothing=1$. We see that the name
``multiplicative game'' is indeed adequate. 

It can be shown by induction that this result holds in general.
\begin{theorem}
Let $v\in\cG(n)$ be a multiplicative game in the above sense. Then its M\"obius
representation satisfies:
\[
m^v(S) = (-1)^{|S|+1}\prod_{i\in [n]\setminus S}\alpha_i.
\]
\end{theorem}

The idea of multiplicative game can be easily generalized to any basis. Let
$\Psi$ be a transform and $(b^\Psi)$ be the corresponding basis. A game
$v\in\cG(n)$ is {\it multiplicative w.r.t. $\Psi$} if it can be
written under the form 
\[
v = \beta(b^\Psi_{\{1\}}-\alpha_1b^\Psi_\varnothing)\times^\Psi(b^\Psi_{\{2\}}- \alpha_2b^\Psi_\varnothing)\times^\Psi\cdots\times^\Psi(b^\Psi_{\{n\}}-\alpha_nb^\Psi_\varnothing)
\]
Then the result of the above theorem remains true and we obtain
\[
\Psi^v(S) = \beta(-1)^{|S|+1}\prod_{i\in[n]\setminus S}\alpha_i.
\]

\section{Cyclotomic games}
Cyclotomic games have a polynomial representation which is a cyclotomic
polynomial. We start by giving some background on cyclotomic polynomials.

\subsection{Cyclotomic polynomials}
Consider for a given $n\in\NN$ the $n$th roots of unity, i.e., the roots of
the polynomial $x^n-1$, which are $e^{2\im\pi\frac{k}{n}}$,
$k=0,1,\ldots,n-1$. The $n$th {\it cyclotomic polynomial} $\Phi_n(x)$ is defined
by
\begin{equation}
\Phi_n(x) = \prod_{\substack{1\leqslant k\leqslant n\\ \gcd(k,n)=1}}(x-e^{2\im\pi\frac{k}{n}}),
\end{equation}
i.e., its roots are the {\it primitive} $n$th roots of unity (those for which
$k$ and $n$ are coprime). Recall that the {\it totient function} of Euler
$\varphi(n)$ gives precisely the number of integers coprime with $n$, from which
we deduce that the degree of $\Phi_n(x)$ is $\varphi(n)$.

We give below the first cyclotomic polynomials:
\begin{align*}
  \Phi_1(x) = & x-1\\
  \Phi_2(x) = & x+1\\
  \Phi_3(x) = & x^2+x+1\\
  \Phi_4(x) = & x^2 +1\\
  \Phi_5(x) = & x^4+x^3+x^2+x+1\\
  \Phi_6(x) = & x^2-x+1\\
  \Phi_7(x) = & x^6+x^5+x^4+x^3+x^2+x+1\\
  \Phi_8(x) = & x^4+1.
\end{align*}

Let us give some properties of cyclotomic polynomials.
\begin{theorem}\label{th:cyclo}
  \begin{enumerate}
  \item $\Phi_n(x)$ is irreducible, monic, with even degree, and has integer coefficients,
    $\forall n\in\NN$.
  \item $\Phi_n(x)$ does not divide $x^k-1$, for every $k<n$, for every $n\in
    \NN$.
  \item $\Phi_n(x)$ is a palindrome, i.e., writing $\Phi_n(x)=a_kx^k+\cdots
    a_1x+a_0$ yields $a_n=a_0$, $a_{n-1}=a_1$, etc.
  \item $\Phi_{2^k}(x)=x^{2^{k-1}}+1$ for every $k\in \NN$.
  \item $\Phi_p(x)=\sum_{k=0}^{p-1}x^k$ when $k$ is a prime number.
  \item $\Phi_{2n}(x)=\Phi_n(-x)$ for $n>1$ and odd.
  \end{enumerate}
\end{theorem}

One can show that moreover, cyclotomic polynomials permit to factorize
$x^n-1$. Indeed,
\begin{align}
  x^n-1 & = \prod_{1\leqslant k\leqslant n}(x - e^{2\im\pi\frac{k}{n}})\nonumber \\
   & = \prod_{d\mid n}\prod_{\substack{1\leqslant k\leqslant\\ \gcd(k,n)=d} }(x
  - e^{2\im\pi\frac{k}{n}})\nonumber\\
  & = \prod_{d\mid n}\Phi_{\frac{n}{d}}(x) = \prod_{d\mid n}\Phi_d(x)\label{eq:cyclo}
\end{align}

\subsection{Cyclotomic games}
A game $v\in \cG$ is said to be {\it cyclotomic} w.r.t. the transform $\Psi$ if
its polynomial representation $P^{v,\Psi}(x)=\Phi_n(x)$ for some $n\in \NN$. From
Theorem~\ref{th:cyclo} (1), we deduce that the canonical representative of
$\Phi_n(x)$ is a game with $\lceil\log_2(\varphi(n))\rceil$
players. We denote by $c_n^\Psi\in\cG$ the canonical representative corresponding to
the $n$th cyclotomic polynomial w.r.t. the transform $\Psi$.

Table~\ref{tab:cyclo} below gives the number of players for the first 30
cyclotomic games. 
\begin{table}[htb]
  \begin{center}
    \begin{tabular}{|c|ccccccccccccccc|}\hline
      $n$th cyclotomic game & 1 & 2 & 3 & 4 & 5 & 6 & 7 & 8 & 9 & 10 & 11 & 12 &
      13 & 14 & 15 \\ \hline
      nb of players & 1 & 1 & 2 & 2 & 3 & 2 & 3 & 3 & 3 & 3 & 4 & 3 & 4 & 3 &
      4\\ \hline\hline
      $n$th cyclotomic game & 16 & 17 & 18 & 19 & 20 & 21 & 22 & 23 & 24 & 25 &
      26 & 27 & 28 & 29 & 30 \\ \hline
      nb of players & 4 & 5 & 3 & 5 & 4 & 4 & 4 & 5 & 4 & 5 & 4 & 5 & 4 & 5 &
      4\\ \hline
    \end{tabular}
  \end{center}
  \caption{Number of players for the 30 first cyclotomic games}
  \label{tab:cyclo}
\end{table}
From the definition and properties in Theorem~\ref{th:cyclo}, we deduce the
following properties for the cyclotomic games.
\begin{theorem}
  The cyclotomic games $c_1^\Psi,c_2^\Psi,\ldots $ have the following properties:
  \begin{enumerate}
  \item $\ga(c_n^\Psi) = \left\{e^{2\im\pi\frac{k}{n}}, \gcd(k,n)=1\right\}$, $\forall n\in \NN$.
  \item $c_n^\Psi$ is irreducible and the last $S$ in the natural order
    s.t. $\Psi^{c_n^\Psi}(S)\neq 0$ does not contain player 1.
  \item $c_{2^k}^\Psi = b^\Psi_{\{k\}}+b^\Psi_{\varnothing}$, $k=1,2,\ldots$.
  \item $c_p^\Psi = \sum_{S\subset \eta^{-1}(p)}b^\Psi_{S}$ for every prime
    number $p$.
  \item \[\Psi^{c_{2n}^\Psi}(S) = \begin{cases}  \Psi^{c_n^\Psi}(S), & \text{ if
    } S\not\ni 1\\
    -\Psi^{c_n^\Psi}(S), & \text{ if } S\ni 1.
    \end{cases}\]
  \end{enumerate}
\end{theorem}

Eq. (\ref{eq:cyclo}) gives the factorization of games whose polynomial
representation is $x^n-1$. Taking the M\"obius basis, this corresponds to the
game $v = u_S - u_{\varnothing}$ for some $S\in\bigcup_{n\in\NN}2^{[n]}$. If the
identity basis is chosen, $v$ is a game defined by $v(S)=1$ for some
$S\in\bigcup_{n\in\NN}2^{[n]}$, $v(\varnothing)=-1$, and $v(T)=0$ otherwise. The
following example illustrates the processus of factorization.
\begin{example}
  Let us consider the M\"obius basis and the game with 4 players
  \[
v = u_{34} - u_{\varnothing}
\]
whose polynomial is $P^{v,m}(x) = x^{12}-1$. Using (\ref{eq:cyclo}) we find
\begin{align*}
  P^{v,m}(x) &= \Phi_1(x)\Phi_2(x)\Phi_3(x)\Phi_4(x)\Phi_6(x)\Phi_{12}(x)\\
  &= (x-1)(x+1)(x^2+x+1)(x^2+1)(x^2-x+1)(x^4-x^2+1).
\end{align*}
The 4-player game has been factored into games of 1, 2 and 3 players, which are:
\begin{align*}
  c_1^m &= u_{1}-u_{\varnothing}\\
  c_2^m &= u_{1}+u_{\varnothing}\\
  c_3^m &= u_{2}+u_{1}+u_{\varnothing}\\
  c_4^m &= u_{2}+u_\varnothing\\
  c_6^m &= u_{2}-u_{1}+u_{\varnothing}\\
  c_{12}^m &= u_{3}-u_{2}+u_{\varnothing}.
\end{align*}
Considering the identity basis or any other works the same way: just replace $u_S$ by
the corresponding basis.

\hfill$\Diamond$

\end{example}

\section{Concluding remarks}
We have proposed a polynomial view of games, where games are represented by an
abstract polynomial, or by its roots as well. This opens a new area of cooperative game theory, where the
classical notions defined for games become irrelevant, forcing the consideration of new
operations on games, like multiplication, and new families, like the
multiplicative games and the cyclotomic games. We believe that these are just
the first steps of an exciting trip into an unknown realm. 

\bibliographystyle{plain}

\bibliography{../BIB/fuzzy,../BIB/grabisch,../BIB/general}

\end{document}